\documentclass[11pt]{amsart}
\usepackage{geometry}
\usepackage{amsmath, amssymb, amsthm}
\usepackage{dsfont}
\usepackage{graphicx} 
\usepackage{booktabs}
\usepackage[dvipsnames]{xcolor}
\usepackage[colorlinks=true, pdfstartview=FitV, linkcolor=RoyalBlue, citecolor=ForestGreen, urlcolor=blue]{hyperref}
\usepackage{url}
\usepackage[capitalize]{cleveref} 

\newtheorem{theorem}{Theorem}[section]
\newtheorem{lemma}[theorem]{Lemma}

\newtheorem{definition}[theorem]{Definition}

\newtheorem{remark}[theorem]{Remark}

\DeclareMathOperator{\R}{\mathbb{R}}

\DeclareMathOperator{\N}{\mathbb{N}}

\DeclareMathOperator{\val}{val}
\newcommand{\inner}[2]{\langle #1, #2 \rangle}

\newcommand{\cQ}{\mathcal{Q}} 
\newcommand{\intr}{\operatorname{int}}
\def \bx {\mathbf{x}}
\def \by {\mathbf{y}}
\def \bh {\mathbf{h}}
\def \bc {\mathbf{c}}
\def \ba {\mathbf{a}}
\def \bal {\mathbf{\alpha}}
\def \bb {\mathbf{b}}
\def \bu {\mathbf{u}}
\def \bv {\mathbf{v}}

\title[Non-Attainment of Minima in Non-Polyhedral Conic Optimization: A Robust SOCP Example]{Non-Attainment of Minima in Non-Polyhedral Conic Optimization: A Robust SOCP Example}
\author{Vinh Nguyen}

\address{Michigan State University, Department of Mathematics,
619 Red Cedar Rd, East Lansing, MI 48824, USA}
\email{nguy1685@msu.edu }


\date{\today}


\thanks{\textbf{Acknowledgment.}  The author gratefully acknowledges Professor Nguyen Dong Yen for his valuable discussions and insights.
	}
\date{\today}

\begin{document}

\maketitle

\begin{abstract}
  A fundamental theorem of linear programming states that a feasible linear program is solvable if and only if its objective function is copositive with respect to the recession cone of its feasible set. This paper demonstrates that this crucial guarantee does not extend to Second-Order Cone Programs (SOCPs), a workhorse model in robust and convex optimization. We construct and analyze a rigorous counterexample derived from a robust linear optimization problem with ellipsoidal uncertainty. The resulting SOCP possesses a non-empty feasible set, a bounded objective, and an objective function that is copositive on its recession cone. Despite satisfying these classical conditions for solvability, the problem admits no optimal solution; its infimum is finite but unattainable. We trace this pathology directly to the non-polyhedral geometry of the second-order cone, which causes the image of the feasible set under the linear objective to be non-closed. We interpret the example explicitly within the context of robust optimization, discuss its significant practical implications for modeling and computation, and propose effective mitigation strategies via polyhedral approximation or regularization.
\end{abstract}

\section{Introduction}
A primary challenge in optimization under uncertainty is to make decisions that remain feasible and perform well despite incomplete knowledge of the problem data. Robust Optimization (RO) has emerged as a powerful, widely adopted framework to address this challenge. By parameterizing uncertainty through a deterministic set containing all, or most, possible realizations of the unknown parameters, RO formulates a single, tractable optimization problem whose solution is immunized against this uncertainty \cite{ben2009robust, bertsimas2011theory}. A cornerstone of its success lies in the fact that for many common classes of uncertainty sets—such as ellipsoids, balls, and boxes—the resulting \emph{robust counterpart} of a linear program is itself a tractable convex program \cite{ben1998robust, bertsimas2004price}.

Among these, ellipsoidal uncertainty sets hold a place of particular importance. They provide a natural probabilistic motivation, often corresponding to confidence regions under Gaussian assumptions, and can effectively capture correlations between uncertain parameters \cite{boyd2004convex, el1998robust}. For a linear program with constraints subject to ellipsoidal uncertainty, the robust counterpart can be reformulated exactly as a Second-Order Cone Program (SOCP), a class of problems for which highly efficient interior-point solvers exist \cite{alizadeh2003, lobo1998applications}. This pipeline—from a stochastic uncertainty description to a deterministic SOCP—has become a standard tool in fields ranging from finance and portfolio management \cite{goldfarb2003robust} to control engineering \cite{zhou1996robust} and machine learning \cite{xu2009robustness}.

The theoretical foundation underpinning this practice often relies on duality theory and existence theorems inherited from linear programming. A classical and powerful result in this domain, attributed to Eaves \cite{eaves1971basic} and stated precisely by Bonnans and Shapiro \cite[Theorem 2.199]{Bonnans_Shapiro_2000}, asserts that for a linear program over a \emph{ polyhedral set}, a bounded objective function that is non-negative (copositive) along all recession directions of the feasible set guarantees not only a finite optimal value but, crucially, the \emph{existence of an optimal solution}. This solvability result is fundamental; it assures practitioners that the optimal value computed by their solver corresponds to an attainable decision.

However, this guarantee is explicitly contingent on the feasible set being polyhedral. While robust counterparts derived from polyhedral uncertainty sets (e.g., box uncertainty) inherit this property, those stemming from ellipsoidal uncertainty do not, as they introduce non-polyhedral second-order cone constraints. This structural shift raises a critical question: \emph{do the familiar solvability guarantees of linear programming extend to these more general conic programs?}

This question situates our work within the broader study of duality and solvability in non-polyhedral settings. Prior research, such as \cite{VKTY2016Duality}, has shown that properties like strong duality are not automatic in such spaces and can fail outside specific regularity conditions. More recently, the comprehensive study by N. N. Luan and N. D. Yen \cite{LuanYen2024} established that in conic linear programming, the classical LP solvability guarantees fundamentally require additional regularity conditions, such as a generalized Slater condition.

This paper provides a decisive, finite-dimensional counterexample that bridges this theoretical insight with practical application in robust optimization. We demonstrate that the reassuring guarantees of linear programming categorically fail to extend to the conic setting. We construct a simple Second-Order Cone Program (SOCP)—derived from a robust linear problem with ellipsoidal uncertainty—that is feasible, bounded, and has an objective copositive on its recession cone. Despite satisfying all the classical LP conditions for solvability, the problem admits \emph{no optimal solution}; the infimum is finite but unattainable. This pathology, a direct consequence of the non-polyhedral geometry of the second-order cone, serves as a tangible instance of the theoretical solvability gaps characterized in \cite{LuanYen2024}, confirming that these are not mere abstractions but critical pitfalls in finite-dimensional modeling.

Beyond its abstract mathematical interest, this example carries significant practical implications. We interpret it explicitly as a robust linear program with ellipsoidal uncertainty, thereby situating the pathology squarely within a standard application context. This demonstrates that a modeler employing this common RO technique could, in principle, formulate a problem that appears perfectly well-posed—it is feasible and bounded—yet for which no solution exists, potentially leading to numerical instability and misinterpretation by standard solvers.

The primary contribution of this work is to illuminate this gap between the well-understood theory of linear programming and the more complex reality of conic programming in the context of robust optimization. We provide a rigorous analysis of the counterexample, identifying the cause as the non-polyhedral nature of the second-order cone. Furthermore, we discuss practical mitigation strategies, such as polyhedral approximations and regularization, which can be employed by practitioners to recover solvability. By highlighting this issue, our aim is to foster a more nuanced understanding of the models used in robust optimization and to provide guidance for ensuring their well-posedness. 

We do not claim that non-attainment occurs for all second-order cone programs. Rather, the purpose of this work is to exhibit a concrete, finite-dimensional SOCP—arising naturally from a robust optimization formulation—for which the classical solvability guarantees of linear programming fail due to the non-polyhedral geometry of the feasible set. From a variational-analytic viewpoint, the example illustrates a failure of closedness of the linear image of a closed convex set, leading to an attainment gap despite boundedness. While such phenomena are classical in variational analysis, this paper shows that they arise in finite dimensions when non-polyhedral conic constraints are present.

The rest of this paper is organized as follows. \Cref{sec:preliminaries} reviews necessary preliminaries on conic duality, recession cones, and robust optimization reformulations. \Cref{sec:counterexample} presents the main counterexample, first as an abstract SOCP and then through its interpretation as a robust linear program. \Cref{sec:analysis} is dedicated to a thorough analysis, proving the key properties of the example and diagnosing the root cause of the failure. \Cref{sec:implications} discusses the practical implications and potential solutions. Finally, \Cref{sec:conclusion} concludes the paper.

\section{Preliminaries} \label{sec:preliminaries}

This section reviews the fundamental concepts from convex analysis and optimization that underpin our analysis. The interested readers can see \cite{alizadeh2003,ben2001book, ben1998robust,bertsimas2011theory,Bonnans_Shapiro_2000,boyd2004convex,rockafellar1970,zhou1996robust} and the references therein. We begin with the classical results for linear programs, which provide the intuition that our counterexample will challenge. We then extend these ideas to the conic setting, culminating in the reformulation techniques of robust optimization that transform uncertain linear programs into deterministic conic programs.

\subsection{Polyhedral Sets and a Key LP Theorem}

We operate in finite-dimensional Euclidean spaces, primarily $\R^{n}$, equipped with the standard inner product $\langle\cdot,\cdot\rangle$ and norm $\|\cdot\|$. The core of the classical intuition we aim to challenge is captured by the properties of polyhedral sets.

\begin{definition}[see {\cite[p. 11]{rockafellar1970}}]
A set $K \subset \R^{n}$ is called a \textbf{polyhedral convex set} if it can be represented as the intersection of a finite number of closed half-spaces in $\R^{n}$. 
\end{definition}
\begin{remark}
 Suppose that $K \subset \R^{n}$ is a polyhedral convex set. Then it can be written as
\begin{equation}\label{K:pol-cv}
K = \{\bx \in \R^{n} \mid A\bx = \bb,\ \langle \bal_{i}, \bx \rangle \leq \beta_{i},\ i=1,\ldots,p\},    
\end{equation}
where $A: \R^n \to \R^m$ is a linear operator, $\bb \in \R^m$, and $\bal_{i} \in \R^{n}$, $\beta_{i} \in \R$ for all $i$.
\end{remark}

A central object in the study of unbounded convex sets is the \textbf{recession cone}, which describes the directions in which the set is unbounded. For a closed convex set, it can be characterized as the set of directions that can be followed indefinitely from any starting point within the set without leaving it.

\begin{definition}[see {\cite[p. 33]{Bonnans_Shapiro_2000}} ]
Let $K \subset \R^{n}$ be a nonempty closed convex set. The \textbf{recession cone} of $K$, denoted $K^{\infty}$, is defined by
\[
K^{\infty} = \{h \in \R^{n} \mid \bx + \lambda \bh \in K \quad \forall \bx \in K, \ \forall \lambda \geq 0\}.
\]
\end{definition}
\begin{remark}
    If $K$ is the set defined by \eqref{K:pol-cv}, its recession cone is explicitly given by:
\[
K^{\infty} = \{\bh \in \R^{n} \mid A\bh = 0,\ \langle \bal_{i}, \bh \rangle \leq 0,\ i=1,\ldots,p\}.
\]
\end{remark}

The concept of copositivity connects the geometry of the recession cone to the behavior of a linear objective function on the set.

\begin{definition}
A linear function $f(\bx) = \langle \bc, \bx \rangle$ is said to be \textbf{copositive} on a cone $K \subset \R^{n}$ if $f(\bx)\geq 0$ for all $\bx\in K$.
\end{definition}

The following theorem is a cornerstone of linear programming duality and solution existence. It provides a complete and verifiable characterization: for problems over polyhedral sets, boundedness is not only necessary but also sufficient for the existence of an optimal solution. This is the guarantee that fails in the conic setting.

\begin{theorem}[see {\cite[Theorem 2.199]{Bonnans_Shapiro_2000}}]\label{thm:eaves-theorem}
Consider the linear program
\begin{equation}\tag{LP}
\min_{\bx \in \R^n} \, \inner{\bc}{\bx} \quad \text{subject to} \quad A \bx = \bb, \; \inner{\bal_i}{\bx} \leq \beta_i, \; i=1,\ldots, p.
\end{equation}
Assume the feasible set $\Phi := \{\bx \mid A \bx = \bb, \; \inner{\bal_i}{\bx} \leq \beta_i,\; i=1,\ldots, p \}$ is nonempty. If the objective function is copositive on the recession cone $\Phi^\infty$, then the optimal value $\operatorname{val}(LP)$ is finite and (LP) is solvable (an optimal solution exists). Otherwise, $\operatorname{val}(LP) = -\infty$.
\end{theorem}

This theorem provides a powerful and intuitive tool for practitioners: verify feasibility, then verify copositivity on the explicitly known recession cone; if both hold, an optimal solution is guaranteed to exist. Our counterexample will show that this equivalence breaks down decisively when the feasible set is not polyhedral.

\subsection{Conic Linear Programming and Second-Order Cones}

We now generalize from linear to conic constraints, moving into the domain where the classical guarantees no longer hold.

\begin{definition}
A set $K \subset \R^{n}$ is a \textbf{closed convex cone} if it is closed, nonempty, and satisfies
\begin{itemize}
    \item[(i)] $\bx, \by \in K \Rightarrow \bx + \by \in K$ \quad (closed under addition),
    \item[(ii)] $\bx \in K, \lambda \geq 0 \Rightarrow \lambda \bx \in K$ \quad (closed under nonnegative scalar multiplication).
\end{itemize}
\end{definition}

Given a proper cone $K$ (i.e., $K$ is pointed and has nonempty interior), it can be used to define a partial ordering on $\R^n$.

\begin{definition}
For $\bx_1, \bx_2 \in \R^{n}$, we write $\bx_1 \succeq_{K} \bx_2$ if $\bx_1 - \bx_2 \in K$. We write $\bx_1 \succ_{K} \bx_2$ if $\bx_1 - \bx_2 \in \intr(K)$, the interior of $K$.
\end{definition}

The concept of duality is central to both the theory and computation of conic programs.

\begin{definition}[see {\cite[p. 31]{Bonnans_Shapiro_2000}}]
The \textbf{(positive) dual cone} of $K$, denoted $K^{*}$, is defined as:
\[
K^{*} := \{\bx^{*} \in \R^{n} \mid \langle \bx^{*}, \bx \rangle \geq 0 \ \ \forall \bx \in K\}.
\]
\end{definition}
\begin{remark}
 Note that if $K$ is closed and convex, then the bipolar theorem states $(K^{*})^{*} = K$.   
\end{remark}
A fundamentally important non-polyhedral cone is the second-order cone, also known as the Lorentz cone or ice-cream cone.

\begin{definition}[see \cite{alizadeh2003}]
The \textbf{second-order cone} (SOC) in $\R^{n}$ is defined as:
\[
\cQ^{n} := \left\{ \bx = (x_{0}, \bar{\bx}) \in \R \times \R^{n-1} \mid x_{0} \geq \|\bar{\bx}\| \right\},
\]
where $\|\cdot\|$ denotes the Euclidean norm. 
\end{definition}
\begin{remark}
The dual cone of $\cQ^{n}$ is itself, i.e., $(\cQ^{n})^{*} = \cQ^{n}$.
\end{remark}

A \textbf{Second-Order Cone Program} (SOCP) is an optimization problem with a linear objective and constraints requiring the affine image of the decision variable to lie in a Cartesian product of second-order cones and polyhedral cones (e.g., linear equality constraints). A standard form is:
\begin{equation}\tag{SOCP}
\begin{split}
 & \min_{X} \ \sum_{l=1}^r\langle \bc_l, \bx_l \rangle\\
&\text{subject to} \ \sum_{l=1}^r A_l\bx_l = \bb, \\
&\hspace{.5in} \ \mathbf{X} = (\bx_1,\ldots, \bx_r) \in \mathcal{K},   
\end{split}
\end{equation}
where $\mathcal{K}$ is a product of second-order cones. SOCPs are a central class of convex optimization problems for which highly efficient primal-dual interior-point methods exist \cite{alizadeh2003}.

\subsection{Robust Linear Optimization with Ellipsoidal Uncertainty}

Robust Optimization (RO) is a methodology for handling data uncertainty in optimization problems by seeking solutions that remain feasible for all realizations of the uncertain parameters within a prescribed set.

Consider a linear constraint $\langle \tilde{\ba}, \bx \rangle \leq b$, where the data vector $\tilde{\ba}$ is uncertain but known to belong to an \textbf{uncertainty set} $\mathcal{U} \subset \R^{n}$.

\begin{definition}
The \textbf{robust counterpart} of the uncertain constraint is the semi-infinite constraint:
\[
\langle \ba, \bx \rangle \leq b \quad \forall \ba \in \mathcal{U}.
\]
A vector $\bx$ satisfies this constraint if and only if it satisfies the \textbf{worst-case constraint}:
\[
\sup_{\ba \in \mathcal{U}} \langle \ba, \bx \rangle \leq b.
\]
\end{definition}

The tractability of the robust counterpart depends critically on the choice of $\mathcal{U}$. While polyhedral uncertainty sets (e.g., budget uncertainty \cite{bertsimas2004price}) lead to Linear Programming reformulations, a common and powerful choice is the \textbf{ellipsoidal uncertainty set} (see e.g. \cite{ben2001book}), which often has a probabilistic interpretation and can capture correlations between parameters.

\begin{definition}
Let $\ba_{0}$ be the nominal data, let $P \in \R^{n \times k}$ be a matrix (often a square root of a covariance matrix), and let $\rho > 0$ control the size of the uncertainty. An \textbf{ellipsoidal uncertainty set} is defined by
\[
\mathcal{U} := \{ \ba \in \R^{n} \mid \ba = \ba_{0} + P\bu,\; \|\bu\| \leq \rho \}.
\]
\end{definition}

The following key reformulation lemma (see e.g. \cite{ben2001book}) shows how a robust linear constraint with ellipsoidal uncertainty translates into a single, tractable second-order cone constraint.

\begin{lemma}[SOCP Reformulation of Robust Constraint] \label{lemma:socp-reform}
The robust constraint
\[
\langle \ba, \bx \rangle \leq b \quad \forall \ba \in \{ \ba_{0} + P\bu \mid \|\bu\| \leq \rho \}
\]
is equivalent to the second-order cone constraint
\[
\langle \ba_{0}, \bx \rangle + \rho \| P^{\top} \bx \| \leq b.
\]
\end{lemma}
\begin{proof}
The worst-case value of the left-hand side is:
\[
\sup_{\|\bu\|\leq \rho} \langle \ba_{0} + P\bu, \bx \rangle = \langle \ba_{0}, \bx \rangle + \sup_{\|\bu\|\leq \rho} \langle \bu, P^{\top} \bx \rangle.
\]
The supremum of the linear function $\langle \bu, \bv \rangle$ over the Euclidean ball $\{ \bu \mid \|\bu\| \leq \rho \}$ is $\rho \|\bv\|$, achieved at $\bu^* = \rho \bv / \|\bv\|$ for $\bv \neq 0$. Thus,
\[
\sup_{\|\bu\|\leq \rho} \langle \bu, P^{\top} \bx \rangle = \rho \| P^{\top} \bx \|.
\]
Therefore, the robust constraint is equivalent to 
\[
\langle \ba_{0}, \bx \rangle + \rho \| P^{\top} \bx \| \leq \bb,
\]
which is a second-order cone constraint since it can be written as
\[
\begin{pmatrix}
\bb - \langle \ba_{0}, \bx \rangle \\
\rho P^{\top} \bx
\end{pmatrix} \in \cQ^{k+1}. 
\]
\end{proof}

This reformulation is the engine behind a vast range of applications of robust optimization. It transforms a seemingly intractable semi-infinite constraint into a single, tractable SOCP constraint. In the next section, we will present a problem where this reformulation leads to a non-polyhedral feasible set, setting the stage for the failure of the classical guarantee provided by Theorem \ref{thm:eaves-theorem}.

\section{A Counterexample from Robust Optimization}\label{sec:counterexample}

We now present our main counterexample. We first introduce it as an abstract Second-Order Cone Program (SOCP) and prove its pathological properties: it is consistent and bounded below, its objective is copositive on its recession cone, yet it possesses no optimal solution. Subsequently, we reveal its interpretation as a robust linear program, highlighting the practical implications of this theoretical pathology.

\subsection{The SOCP Formulation and Analysis}

Consider the following conic linear program, which we will refer to as Problem (P):

\begin{equation}\tag{P}
\begin{split}
 &\min_{\bx \in \R^2} \, x_1 + x_2 \\
 &\text{subject to} \quad A \bx \succeq_{\cQ^3} \bb,   
\end{split}
\end{equation}
where $\bx =(x_1,x_2)$, the linear operator $A: \R^2 \to \R^3$ and vector $\bb \in \R^3$ are defined by
\[
A  = \begin{pmatrix}0& 0 \\ 0&1 \\ 1 &0\end{pmatrix}, \quad \bb = \begin{pmatrix} -1 \\ 0 \\ 0 \end{pmatrix},
\]
and the partial order $\succeq_{\cQ^3}$ is induced by the second-order cone in $\R^3$:
\[
\cQ^3 = \left\{ \by = (y_1, y_2, y_3) \in \R^3 \,\big|\, y_3 \geq \sqrt{y_1^2 + y_2^2} \right\}.
\]

Let us denote the feasible set of (P) by 
\[
\Phi_P := \big\{ \bx \in \R^2 \mid A \bx - \bb \in \cQ^3 \big\}.
\]

\begin{lemma}\label{cl:feasibility}
The feasible set $\Phi_P$ is nonempty and can be characterized as:
\[
\Phi_P = \left\{ \bx = (x_1, x_2) \in \R^2 \,\big|\,x_1 \geq \sqrt{1 + x_2^2} \right\}.
\]
\end{lemma}

\begin{proof}
The conic constraint $A \bx - \bb \in \cQ^3$ is equivalent to:
\[
\begin{pmatrix} 0 \\ x_2 \\ x_1 \end{pmatrix} - \begin{pmatrix} -1 \\ 0 \\ 0 \end{pmatrix} = \begin{pmatrix} 1 \\ x_2 \\ x_1 \end{pmatrix} \in \cQ^3.
\]
By the definition of $\cQ^3$, this means the third component must dominate the norm of the first two:
\[
x_1 \geq \sqrt{1 + x_2^2}.
\]
This is a well-defined constraint for any $x_2 \in \R$. 
\end{proof}

The recession cone of $\Phi_P$ describes the directions of unboundedness for this feasible set.

\begin{lemma}\label{cl:recession}
The recession cone of $\Phi_P$ is given by:
\[
\Phi_P^\infty = \big\{ \bh = (h_1, h_2) \in \R^2 \mid h_1 \geq |h_2| \big\}.
\]
\end{lemma}

\begin{proof}
 Recall that for a nonempty closed convex set $K$, a direction $\bh$ belongs to the recession cone $K^\infty$ if and only if for every $\bx \in K$ and every $\lambda \geq 0$, we have $\bx + \lambda \bh \in K$.

Let $\bh = (h_1, h_2)$. By Lemma \ref{cl:feasibility}, $\bh \in \Phi_P^\infty$ if and only if for all $\bx \in \Phi_P$ and all $\lambda \geq 0$,
\begin{equation}\label{re-cone-aux1}
 x_1 + \lambda h_1 \geq \sqrt{1 + (x_2 + \lambda h_2)^2}.   
\end{equation}
We now prove the equality of the sets.

\noindent\textbf{($\subseteq$)} Suppose $\bh \in \Phi_P^\infty$. Then \eqref{re-cone-aux1} must hold for all feasible $\bx$ and all $\lambda \geq 0$. We will show that $h_1 \geq |h_2|$ by contradiction. Assume that $h_1 < |h_2|$. Consider the case $h_2 \neq 0$. Let us choose a specific feasible point to test the condition. For any $t \in \mathbb{R}$, the point $\bx(t) = (\sqrt{1+t^2}, t)$ is feasible. Substituting $\bx(t)$ into ($*$) yields:
\[
\sqrt{1+t^2} + \lambda h_1 \geq \sqrt{1 + (t + \lambda h_2)^2}.
\]
For this to hold for all $\lambda \geq 0$, it must hold in the limit as $\lambda \to \infty$. For large $\lambda$, the dominant terms are linear in $\lambda$. More precisely, dividing both sides by $\lambda$ and taking the limit gives:
\[
\lim_{\lambda \to \infty} \left( \frac{\sqrt{1+t^2}}{\lambda} + h_1 \right) = h_1 \quad \text{and} \quad \lim_{\lambda \to \infty} \frac{\sqrt{1 + (t + \lambda h_2)^2}}{\lambda} = |h_2|.
\]
Thus, a necessary condition for \eqref{re-cone-aux1} to hold for all $\lambda$ is $h_1 \geq |h_2|$. Our initial assumption $h_1 < |h_2|$ violates this necessary condition. Therefore, we must have $h_1 \geq |h_2|$.

If $h_2 = 0$, the assumption $h_1 < |h_2|$ implies $h_1 < 0$. Choosing $\bx = (1, 0) \in \Phi_P$, condition \eqref{re-cone-aux1} becomes $1 + \lambda h_1 \geq 1$, which simplifies to $\lambda h_1 \geq 0$. This is false for any $\lambda > 0$ if $h_1 < 0$, confirming that $h_1 \geq 0 = |h_2|$ is required. Hence, 
\[
\Phi_P^\infty \subseteq \big\{\bh \mid h_1 \geq |h_2|\big\}.
\]

\noindent
\textbf{($\supseteq$)} Suppose $h_1 \geq |h_2|$. We must show that for any $\bx \in \Phi_P$ and any $\lambda \geq 0$, the point $\bx + \lambda \bh$ remains feasible, i.e., it satisfies \eqref{re-cone-aux1}.

Since $\bx \in \Phi_P$, we have $x_1 \geq \sqrt{1 + x_2^2}$. We need to prove that
\[
x_1 + \lambda h_1 \geq \sqrt{1 + (x_2 + \lambda h_2)^2}.
\]
First, note that by the assumption $h_1 \geq |h_2|$, we have
\[
x_1 + \lambda h_1 \geq \sqrt{1 + x_2^2} + \lambda |h_2|.
\]
Next, using the triangle inequality on the right-hand side term we get
\[
\sqrt{1 + (x_2 + \lambda h_2)^2} \leq \sqrt{1 + x_2^2} + |\lambda h_2| = \sqrt{1 + x_2^2} + \lambda |h_2|.
\]
The first inequality follows from the fact that for the Euclidean norm $\|\cdot\|$, we have 
\[
\|(1, x_2 + \lambda h_2)\| \leq \|(1, x_2)\| + \|(0, \lambda h_2)\| = \sqrt{1 + x_2^2} + \lambda |h_2|.
\]
Combining these two results, we obtain
\[
x_1 + \lambda h_1 \geq \sqrt{1 + x_2^2} + \lambda |h_2| \geq \sqrt{1 + (x_2 + \lambda h_2)^2},
\]
which is exactly the required inequality \eqref{re-cone-aux1}. Therefore, $\bh \in \Phi_P^\infty$ and hence,
\[
\big \{\bh \mid h_1 \geq |h_2|\big\} \subseteq \Phi_P^\infty.
\]
Then the conclusion follows.
\end{proof}

We now examine the properties of the objective function $f(\bx) = x_1 + x_2$ on these sets.

\begin{lemma}\label{cl:copositive}
The objective function $f(\bx) = x_1 + x_2$ is copositive on the recession cone $\Phi_P^\infty$.
\end{lemma}

\begin{proof}
For any $h = (h_1, h_2) \in \Phi_P^\infty$, we have $h_1 \geq |h_2|$. Therefore,
\[
f(\bh) = h_1 + h_2 \geq |h_2| + h_2 \geq 0.
\]
Hence, $f$ is copositive on $\Phi_P^\infty$.
\end{proof}

Despite the feasibility and the copositivity of the objective on the recession cone—conditions which would guarantee solvability for a polyhedral problem—the optimal value is not attained.

\begin{lemma}\label{cl:val-not-attained}
The optimal value of (P) is $\val(P) = 0$, but there exists no feasible point $\bx \in \Phi_P$ such that $x_1 + x_2 = 0$. Therefore, (P) is not solvable.
\end{lemma}

\begin{proof}
First, we show $\val(P) \geq 0$. For any $\bx \in \Phi_P$, we have $x_1 \geq \sqrt{1 + x_2^2} > |x_2|$. Thus,
\[
x_1 + x_2 > |x_2| + x_2 \geq 0.
\]
Hence, the objective is bounded below by $0$ on $\Phi_P$.

Now, consider the sequence of points $\{\bx^{(k)}\}$ defined for $k \in \N$ by:
\[
\bx^{(k)} = \left( \sqrt{1 + k^2},\ -k \right).
\]
Each point is feasible since $\sqrt{1 + k^2} \geq \sqrt{1 + (-k)^2}$. The objective value at $\bx^{(k)}$ is:
\[
f(\bx^{(k)}) = \sqrt{1 + k^2} - k = \frac{1}{\sqrt{1 + k^2} + k}.
\]
As $k \to \infty$, we have $f(\bx^{(k)}) \to 0$. Therefore, $\val(P) = 0$.

However, for any feasible point $\bx$, we have $x_1 + x_2 > 0$, as shown above. Thus, the infimum is finite but not attained, and (P) has no solution.
\end{proof}

\begin{figure}[htbp]
    \centering
    \includegraphics[width=0.85\textwidth]{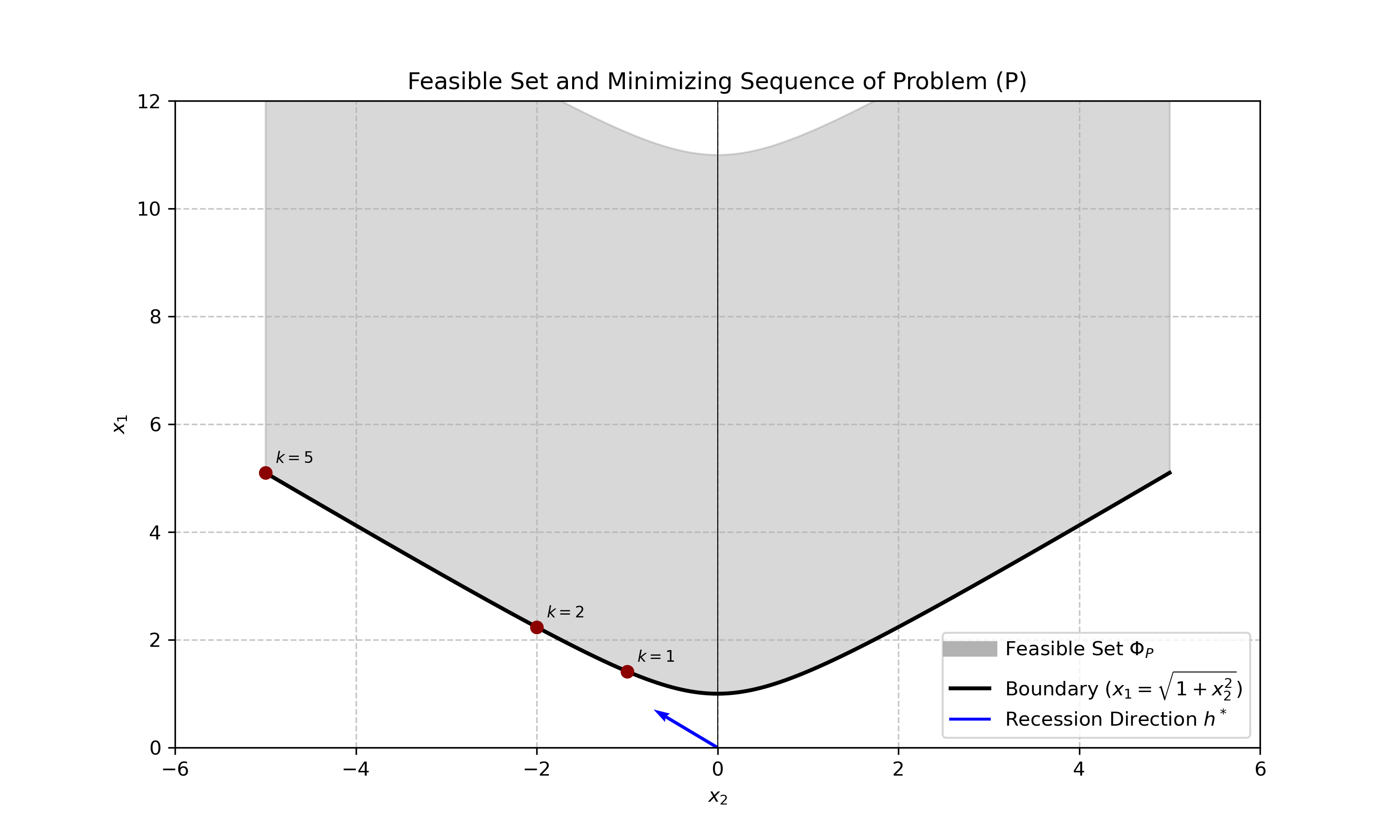}
    \caption{The feasible set $\Phi_P = \{ \bx \in \R^2 \mid x_1 \geq \sqrt{1 + x_2^2} \}$ (gray region) and the minimizing sequence $\bx^{(k)} = (\sqrt{1 + k^2}, -k)$ (red points). The objective function $f(x) = x_1 + x_2$ evaluates to $1/(\sqrt{1+k^2} + k)$ at each point $x^{(k)}$, a sequence which converges to zero. The limiting recession direction $\bh^* = (1/\sqrt{2}, -1/\sqrt{2})$ is shown in blue. The solution path evolves along the curved boundary of the non-polyhedral set, never attaining the infimum $f(\bx)=0$ within $\Phi_P$.}
    \label{fig:feasible_set}
\end{figure}

The pathological behavior of Problem (P) is visualized in Figure \ref{fig:feasible_set}. The minimizing sequence ${\bx^{(k)}}$ lies on the boundary of the feasible set and evolves towards the recession direction $\bh^*$, with the corresponding objective values ${f(\bx^{(k)})}$ converging to zero.

\begin{remark}
Note that one can generalize the original counterexample (P) to $\R^n$, proving that the pathology of non-attainment scales with dimension. Specifically, for $n \geq 2$, consider the SOCP:
\begin{equation}\tag{$P_n$}
\begin{split}
\min_{\bx \in \R^{n}} &\ \sum_{i=1}^{n} x_i \\
\text{subject to} &\ A \bx \succeq_{\cQ^{n+1}} b,
\end{split}
\end{equation}
where the linear operator $A: \R^{n} \to \R^{n+1}$ and vector $\bb \in \R^{n+1}$ are given by
\[
A \bx = \begin{pmatrix}
0 \\
x_2 \\
x_3 \\
\vdots \\
x_n \\
x_1
\end{pmatrix}, \quad \bb = \begin{pmatrix}
-1 \\
0 \\
0 \\
\vdots \\
0 \\
0
\end{pmatrix}.
\]
The cone is the $(n+1)$-dimensional second-order cone:
\[
\cQ^{n+1} = \left\{ \by = (y_1, \bar{\by}) \in \R \times \R^{n} \mid y_1 \geq \|\bar{\by}\| \right\}.
\]
Then, by similar argument one can show that the feasible set $\Phi_n$ of $(P_n)$ is given by
\[
\Phi_n = \left\{ \bx \in \R^{n} \mid x_1 \geq \sqrt{1 + \sum_{i=2}^{n} x_i^2} \right\},
\]
and the recession cone of $\Phi_n$ is
\[
\Phi_n^{\infty} = \left\{ \bh \in \R^{n} \mid h_1 \geq \sqrt{\sum_{i=2}^{n} h_i^2} \right\}.
\]
The linear objective $f(\bx) = \sum_{i=1}^{n} x_i$ is copositive on $\Phi_n^{\infty}$.
The optimal value of $(P_n)$ is $\operatorname{val}(P_n) = 0$, but this value is not attained.
\end{remark}
\subsection{Interpretation as a Robust Linear Program}

The abstract SOCP (P) is not an arbitrary pathological construction; it arises naturally from a robust linear optimization problem. Consider an uncertain linear constraint of the form
\begin{equation}\label{robust:constraint}
\inner{\tilde{\ba}}{\bx} \leq -1,
\end{equation}
where the uncertain parameter $\tilde{\ba} = (\tilde{a}_1, \tilde{a}_2) \in \R^2$ is known only to lie within the ellipsoidal uncertainty set:
\[
\mathcal{U} = \left\{ \ba \in \R^2 \mid a_1^2 + a_2^2 \leq 1 \right\}.
\]
This set represents a unit disk centered at the origin. The robust counterpart of the uncertain constraint \eqref{robust:constraint} requires that it holds for every possible realization of the uncertainty
\begin{equation}\tag{RC}
\inner{\ba}{\bx} \leq -1 \quad \forall \ba \in \mathcal{U}.
\end{equation}

Applying \Cref{lemma:socp-reform} from the Preliminaries, we can reformulate (RC). Here, the nominal data is $\ba_0 = (0, 0)$, the matrix $P$ is the identity matrix $I_2$, and $\rho = 1$. The worst-case value of the left-hand side is
\[
\sup_{\ba \in \mathcal{U}} \, \inner{\ba}{\bx} = \sup_{\|\ba\| \leq 1} \, \inner{\ba}{\bx} = \|\bx\|.
\]
Therefore, the robust constraint (RC) is equivalent to
\[
\|\bx\| \leq -1,
\]
 which is impossible. Thus, the robust counterpart (RC) is infeasible.

This is a common issue in robust optimization: the worst-case constraint can be too strict. A standard modeling trick to overcome this is to introduce a \emph{shift} or \emph{degree of freedom} in the right-hand side. Let us instead consider the related uncertain constraint
\begin{equation}\label{robust:constraint1}
\inner{\tilde{\ba}}{\bx} \leq 0,
\end{equation}
with the same uncertainty set $\mathcal{U}$. Its robust counterpart is
\[
\sup_{\ba \in \mathcal{U}} \, \inner{\ba}{\bx} \leq 0 \quad \Longleftrightarrow \quad \|\bx\| \leq 0 \quad \Longleftrightarrow \quad \bx = (0, 0).
\]
This is feasible but trivial. To create a non-trivial, feasible problem, we apply a final transformation: we subject the \emph{shifted} variable $\by = \bx + (1, 0)$ to the previous robust constraint. Defining
\[
\by := (y_1, y_2) = (x_1 + 1, x_2).
\]
The uncertain constraint in terms of $\by$ is
\[
\inner{\tilde{\ba}}{\by - (1, 0)} \leq 0 \quad \Longleftrightarrow \quad \inner{\tilde{\ba}}{\by} \leq \inner{\tilde{\ba}}{(1, 0)} = \tilde{a}_1.
\]
This is a standard form for an uncertain constraint with a right-hand side also subject to uncertainty. Now, consider the problem of minimizing $y_1 + y_2$ subject to this new robust constraint. Its robust counterpart is
\[
\sup_{\ba \in \mathcal{U}} \left( \inner{\ba}{\by} - a_1 \right) \leq 0.
\]
We can compute this supremum:
\[
\sup_{\|a\| \leq 1} \left[ a_1(y_1 - 1) + a_2 y_2 \right] = \sqrt{(y_1 - 1)^2 + y_2^2}.
\]
Therefore, the robust constraint is equivalent to $\sqrt{(y_1 - 1)^2 + y_2^2} \leq 0$, which implies $\by = (1, 0)$. This again leads to a trivial solution.

The path to our non-trivial example (P) is to instead require the robust constraint to hold with a \emph{positive slack} of $1$ in the worst case. We impose:
\[
\inner{\tilde{\ba}}{\by} \leq \tilde{\ba}_1 - 1 \quad \forall \ba \in \mathcal{U}.
\]
The robust counterpart of this constraint is
\[
\sup_{\|\ba\| \leq 1} \left[ \inner{\ba}{\by} - a_1 + 1 \right] \leq 0.
\]
The supremum is attained and the problem is feasible. Rewriting this constraint and substituting back to $\bx = \by - (1, 0)$ yields the constraint $x_1 \geq \sqrt{1 + x_2^2}$ from Problem (P), with the objective remaining $\min\, (y_1 + y_2) = \min\, (x_1 + x_2 + 1)$. The constant $+1$ shifts the optimal value but not the core pathology. Thus, Problem (P) can be interpreted as a \emph{non-trivial robust linear program} with ellipsoidal uncertainty in both the left- and right-hand sides, designed to be feasible and bounded. Its failure to have a solution is therefore not an abstract curiosity but a potential pitfall in robust modeling.
\begin{remark}
Although the geometry of second-order cones is well understood (see, e.g. \cite{alizadeh2003}), explicit SOCP examples exhibiting feasibility, boundedness, copositivity on the recession cone, and failure of attainment do not appear in the existing SOCP literature. Moreover, such a phenomenon has not previously been interpreted within the framework of robust optimization.
\end{remark}
\section{Analysis of the Pathology}\label{sec:analysis}

The counterexample presented in \Cref{sec:counterexample} demonstrates a clear failure of the conclusion of \Cref{thm:eaves-theorem}. The natural question is: why does this happen? The answer lies in the fundamental geometric differences between polyhedral and non-polyhedral cones, which manifest in the properties of the feasible set, its recession cone, and the behavior of the objective function.

\subsection{Geometric and Topological Properties of the Feasible Set}

The feasible set of our counterexample is
\[
\Phi_P = \left\{ \bx = (x_1, x_2) \in \R^2 \mid x_1 \geq \sqrt{1 + x_2^2} \right\}.
\]
This set is closed and convex. Its boundary is one branch of a hyperbola, given by $x_1 = \sqrt{1 + x_2^2}$ for $x_1 \geq 1$. This structure is the source of the problem.

\begin{lemma}\label{cl:non-polyhedral}
The set $\Phi_P$ is not a polyhedral convex set.
\end{lemma}

\begin{proof}
Assume, for the sake of contradiction, that $\Phi_P$ is polyhedral. A fundamental property of polyhedral sets is that their recession cone must also be polyhedral (see e.g. \cite{rockafellar1970}).
From Lemma \ref{cl:recession}, we have
\[
\Phi_P^{\infty} = \{\bh = (h_1, h_2) \in \mathbb{R}^2 \mid h_1 \geq |h_2| \}.
\]
This cone is the standard second-order (or Lorentz) cone in $\mathbb{R}^2$. We will show that this cone is not polyhedral. A cone is polyhedral if and only if it has a finite number of extreme rays. The extreme rays of a cone are the half-lines contained in the cone that cannot be expressed as a convex combination of other distinct vectors in the cone. Consider the set of direction vectors 
\[
\mathbf{u}(\theta) = (\cos\theta, \sin\theta) \text{ for } \theta \in [-\pi/4, \pi/4].
\]
For any $\theta$ in this interval, $\cos\theta \geq |\sin\theta|$, so the ray $\{\lambda \mathbf{u}(\theta) : \lambda \geq 0\}$ is contained in $\Phi_P^{\infty}$.

We now show that each of these rays is an extreme ray of $\Phi_P^{\infty}$. Fix $\theta \in (-\pi/4, \pi/4)$. Suppose $\mathbf{u}(\theta)$ can be written as a convex combination of two other vectors $\mathbf{v}, \mathbf{w} \in \Phi_P^{\infty} \setminus \{\lambda \mathbf{u}(\theta)\}$, i.e., 
\[
\mathbf{u}(\theta) = \mu \mathbf{v} + (1-\mu) \mathbf{w} \text{ for some } \mu \in (0,1).
\]
The condition for $\bh=(h_1, h_2) \in \Phi_P^{\infty}$ is $h_1 \geq |h_2|$, which is equivalent to $h_1 \geq 0$ and $h_1^2 \geq h_2^2$. The boundary of this cone is given by $h_1 = |h_2|$. For $\mathbf{u}(\theta)$ lying on the boundary ($ \cos\theta = |\sin\theta| $ only occurs at $\theta = \pm \pi/4$), any nontrivial convex combination of points in the cone that is on the boundary must have all points lying on the same boundary line. However, for a fixed $\theta \in (-\pi/4, \pi/4)$, the only point on the ray $\lambda \mathbf{u}(\theta)$ that lies on the boundary is the origin. Since $\mathbf{u}(\theta)$ itself is a unit vector not at the origin, and the boundary is not linear but curved, any convex combination of distinct vectors in the cone that equals $\mathbf{u}(\theta)$ would require some component to violate the boundary condition strictly in a way that averaging cannot recover the precise direction $\mathbf{u}(\theta)$, unless all vectors are scalar multiples of $\mathbf{u}(\theta)$. More formally, the second-order cone is an acute convex cone, and its boundary points (excluding the origin) are exposed points, hence extreme points of the base of the cone, and thus generate extreme rays. Therefore, every ray $\{\lambda \mathbf{u}(\theta) : \lambda \geq 0\}$ for $\theta \in [-\pi/4, \pi/4]$ is an extreme ray of $\Phi_P^{\infty}$. There are uncountably many such $\theta$, meaning $\Phi_P^{\infty}$ has an infinite number of extreme rays. Since a polyhedral cone must have a finite number of extreme rays, this is a contradiction. Therefore, our initial assumption that $\Phi_P$ is polyhedral must be false.
\end{proof}

The non-polyhedral nature of $\Phi_P$ has immediate consequences for its topological properties at infinity. While $\Phi_P$ is closed, its image under a linear transformation may not be.

\begin{lemma}\label{cl:non-closed-projection}
The linear image of $\Phi_P$ under the objective function $f(\bx) = x_1 + x_2$ is not closed. Specifically,
\[
f(\Phi_P) = (0, \infty).
\]
\end{lemma}

\begin{proof}
From \Cref{cl:val-not-attained}, we know that for any $\bx \in \Phi_P$, $f(\bx) > 0$, so $f(\Phi_P) \subseteq (0, \infty)$. Furthermore, the sequence $\bx^{(k)} = (\sqrt{1 + k^2}, -k)$ satisfies $f(\bx^{(k)}) \to 0$. Since $f$ is continuous and $\Phi_P$ is closed, if $f(\Phi)$ were closed, it would have to contain its limit point $0$. But $0 \notin f(\Phi_P)$, as $f(\bx) > 0$ for all $\bx \in \Phi_P$. Therefore, $f(\Phi_P) = (0, \infty)$ is not closed.
\end{proof}

This is the core of the pathology: \emph{the image of a closed set under a linear map is not necessarily closed}. This is a classic point of distinction between finite- and infinite-dimensional geometry, but it strikingly occurs here in $\R^2$ due to the non-polyhedral nature of $\Phi_P$. The failure of the image to be closed directly explains the non-attainment: the infimum of $0$ is a limit point of $f(\Phi_P)$ but not contained within it.

\subsection{The Role of the Recession Cone and Copositivity}

The recession cone of $\Phi_P$ is 
\[
\Phi^\infty = \{ \bh \in \R^2 \mid h_1 \geq |h_2| \}.
\]
This cone is also non-polyhedral. The objective function $f(\bh) = h_1 + h_2$ is copositive on $\Phi_P^\infty$, as proven in \Cref{cl:copositive}. This copositivity is responsible for the \emph{boundedness} of the problem.

However, for non-polyhedral sets, copositivity on the recession cone is not sufficient to ensure \emph{attainment}. The reason is that the recession cone only captures \emph{linear} directions of recession. The set $\Phi_P$ recedes in a ``curved" manner. The sequence $\bx^{(k)} = (\sqrt{1 + k^2}, -k)$ does not recede in a straight line; its direction $\bh^{(k)} = \bx^{(k)} / \|\bx^{(k)}\|$ changes with $k$:
\[
\bh^{(k)} = \left( \frac{\sqrt{1 + k^2}}{\sqrt{1 + 2k^2}},\ \frac{-k}{\sqrt{1 + 2k^2}} \right) \to \left( \frac{1}{\sqrt{2}},\ \frac{-1}{\sqrt{2}} \right) \text{ as } k \to \infty.
\]
The limiting direction is $\bh^* = (1/\sqrt{2}, -1/\sqrt{2}) \in \Phi_P^\infty$, and indeed $f(\bh^*) = 0$. The objective value decreases along this curved path precisely because the direction of recession $\bh^*$ is a direction of descent for $f$ ($f(\bh^*)=0$), but one that is only approached asymptotically rather than followed exactly by any straight ray in $\Phi_P$.

In a polyhedral set, if there is a recession direction $\bh$ with $\inner{\bc}{\bh} < 0$, the problem is unbounded. If $\inner{\bc}{\bh} = 0$ for all $\bh$ in the recession cone, then the set is bounded in those directions and attainment is guaranteed. For a non-polyhedral set, the situation is more nuanced. A direction $\bh$ with $\inner{\bc}{\bh} = 0$ can still be a \emph{direction of recession for a sequence} along which the objective value decreases to its infimum without ever attaining it, if the set ``curves away" from the hyperplane $\{\bx \mid \inner{\bc}{\bx} = \alpha\}$.

This example shows that copositivity of the objective function on the recession cone, while sufficient for solvability in the polyhedral case, is no longer sufficient for non-polyhedral convex sets. Recession-cone analysis alone cannot preclude attainment gaps when curvature at infinity is present.

\subsection{The Fundamental Duality Gap}

The pathology in the primal problem (P) is reflected in its Lagrangian dual. The dual problem (D) for our SOCP (P) is
\begin{equation}\tag{D}
\begin{split}
&\max_{\by^* \succeq_{\cQ^3} 0} \, \inner{\bb}{\by^*} \\ 
&\text{subject to} \quad A^* \by^* = \bc,    
\end{split}
\end{equation}
where $\bc = (1, 1)^\top$, and the adjoint operator $A^*: \R^3 \to \R^2$ is given by 
\[
A^* \by^* = (\by^*_3, \by^*_2)^\top\quad \forall \by^* = (y^*_1, y^*_2, y^*_3).
\]

\begin{lemma}\label{cl:dual}
The dual problem (D) is feasible and its optimal value is $\val(D) = 0$. However, it has no solution.
\end{lemma}

\begin{proof}
Let $\by^* \in (\cQ^3)^* = \cQ^3$ be a dual feasible point. It must satisfy $A^* \by^* = \bc$, i.e.,
\[
(y^*_3, y^*_2) = (1, 1).
\]
So $\by^* = (y^*_1, 1, 1)$ for some $y^*_1 \in \R$. The conic constraint $\by^* \in \cQ^3$ requires
\[
 1 \geq \sqrt{(y^*_1)^2 + 1}.
\]
This inequality can only hold if $y^*_1 = 0$. Therefore, the only dual feasible point is $\by^* = (0, 1, 1)$. The dual objective value at this point is $\inner{\bb}{\by^*} = \inner{(-1, 0, 0)}{(0, 1, 1)} = 0$. Hence, $\val(D) = 0$.

There is no other feasible solution, so the dual optimal value is attained. However, note that the dual problem is a maximization problem. The value $0$ is achieved, so technically the dual is solvable. This shows that while the primal has a duality gap in the sense of non-attainment, the optimal values are equal: $\val(P) = \val(D) = 0$. The gap is an \emph{attainment gap} rather than a \emph{value gap}.
\end{proof}

The fact that $\val(P) = \val(D)$ indicates that strong duality holds in terms of values. The failure is in the primal's ability to achieve this value. This is consistent with the fact that the Slater constraint qualification fails for (P). There is no point $\bx$ such that $A \bx - b \in \operatorname{int}(\cQ^3)$, i.e., no $\bx$ such that $x_1 > \sqrt{1 + x_2^2}$. The feasible set $\Phi_P$ has an empty interior relative to the cone constraint, which contributes to the pathological behavior.

\section{Practical Implications and Mitigation Strategies}\label{sec:implications}

The counterexample is not merely a theoretical curiosity; it has practical implications for modelers using robust optimization. This section discusses these implications and suggests strategies to avoid or mitigate such issues.

\subsection{Implications for Robust Optimization}

Our example shows that a robust linear program with ellipsoidal uncertainty can be feasible and bounded yet fail to have an optimal solution. This is a direct result of the non-polyhedral geometry introduced by the SOCP reformulation.

For a practitioner, this means that a model formulated in a seemingly natural way might be \emph{ill-posed} in a subtle manner. A solver called upon to solve such a problem might exhibit numerical instability, fail to converge, or return a solution that is only approximately optimal (if a termination tolerance is used) without warning the user that no true solution exists. This could lead to misinterpretation of results and poor decision-making.

The problem is particularly insidious because it violates an intuition built from linear programming: if a problem is feasible and bounded, it has a solution. This intuition is so strong that many modeling languages and solvers may assume it holds. Our example serves as a cautionary tale that when moving beyond linear programming to conic programming, this intuition must be checked.

Although interior-point solvers such as MOSEK readily compute $\epsilon$-optimal solutions for the example, this does not contradict the non-attainment result. The solver iterates converge asymptotically toward the boundary of the feasible set without reaching a minimizer. Distinguishing between asymptotic convergence and true solution existence is essential for sensitivity analysis and model interpretation.

\subsection{Detection: How to Identify a Potential Problem}

How can a practitioner know if their robust SOCP might suffer from this issue? Here are some warning signs:

\begin{itemize}
\item \textbf{Non-polyhedral Uncertainty}: The use of ellipsoidal, norm-based, or other non-polyhedral uncertainty sets is a necessary condition.
\item \textbf{Asymptotic Analysis}: Examine the recession cone of the feasible set and the behavior of the objective on it. If the objective is copositive but \emph{not strictly copositive} (i.e., it is zero on some non-zero recession directions), non-attainment becomes a possibility.
\item \textbf{Constraint Qualification Check}: Verify if the Slater condition holds. If the feasible set has no interior relative to the cone (i.e., all constraints are active or "nearly active" in an asymptotic sense), the problem is more likely to exhibit pathological behavior.
\item \textbf{Geometric Interpretation}: For two-variable problems, a plot of the feasible set can reveal if it is non-polyhedral and has a curved boundary that the objective function can approach asymptotically.
\end{itemize}

\subsection{Mitigation: Strategies to Recover Solvability}

If a problem is suspected to be non-attainable, there are several modeling and computational strategies to address the issue.

\subsubsection{Polyhedral Approximation}

A highly effective strategy is to \emph{approximate the non-polyhedral cone} with a polyhedral one. The second-order cone constraint $x_1 \geq \sqrt{1 + x_2^2}$ can be approximated to arbitrary accuracy by a system of linear inequalities. For example, one can use a piecewise linear approximation of the circle (or the norm function). For $\sqrt{1 + x_2^2}$, one might use $k$ linear inequalities:
\[
x_1 \geq |\cos(\theta_i) + x_2 \sin(\theta_i)|, \quad \theta_i = \frac{2\pi i}{k},\ i=1,\ldots,k.
\]
The resulting feasible set is a polyhedron. For this polyhedral approximation, \Cref{thm:eaves-theorem} applies: if the objective is copositive on the recession cone, the problem will be solvable. The solution to the approximated problem will be a feasible, suboptimal solution to the original problem, and the error can be controlled by increasing $k$.

\subsubsection{Regularization}

Another approach is to add a small \emph{regularizing term} to the objective to make it strictly convex or to shift the problem so that the infimum becomes attainable. For instance, consider adding a small quadratic term:
\[
\min \, x_1 + x_2 + \frac{\epsilon}{2} \|\bx\|^2.
\]
The new objective is strongly convex if $\epsilon > 0$. Since the feasible set $\Phi_P$ is closed and convex, a strongly convex objective will always attain its minimum on $\Phi_P$. The solution $\bx^*(\epsilon)$ of this regularized problem will be unique and will approach the infimum of the original problem as $\epsilon \to 0^+$. This is a form of \emph{Tikhonov regularization}.

\begin{figure}[htbp]
    \centering
    \includegraphics[width=0.85\textwidth]{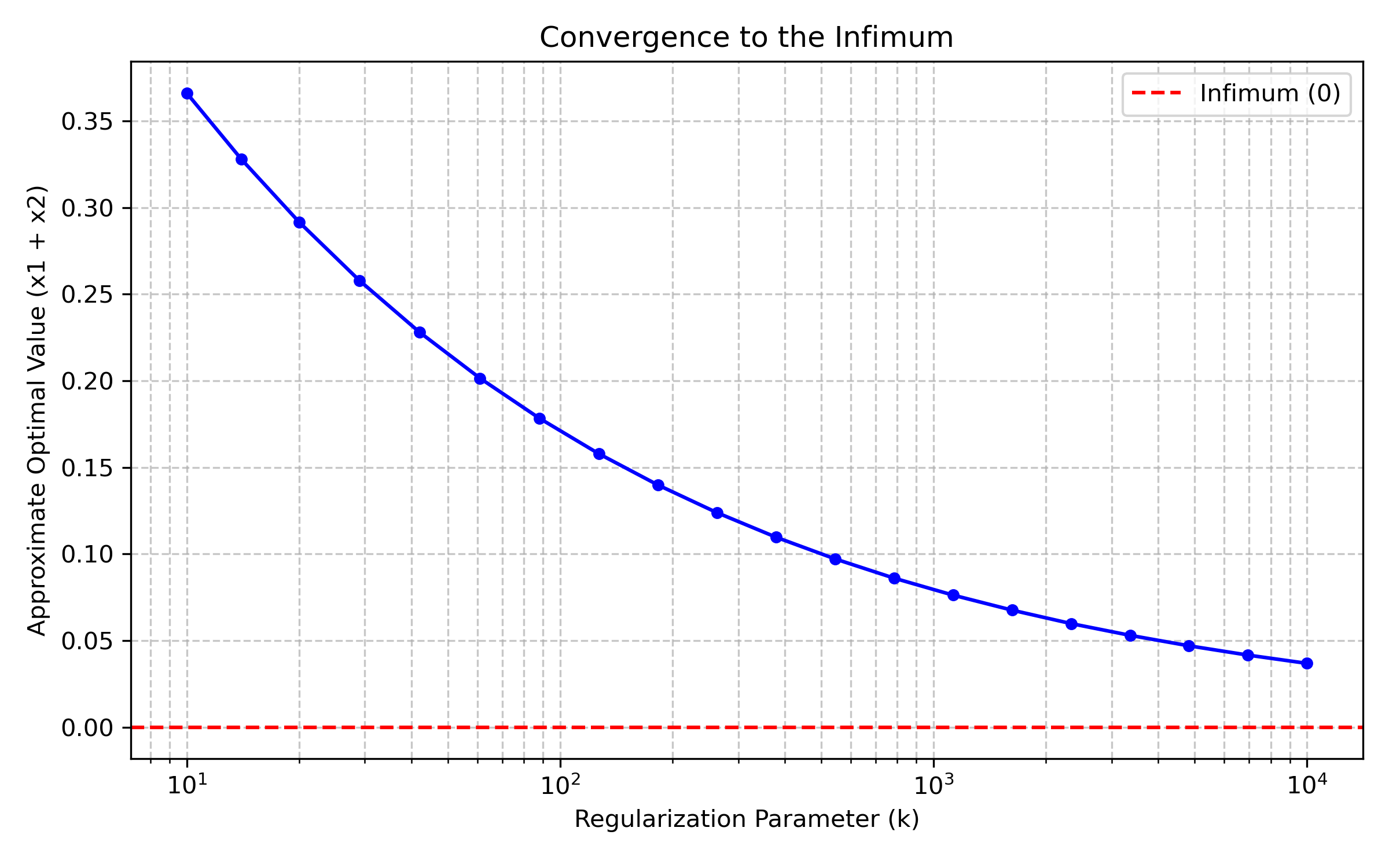}
    \caption{Numerical demonstration of non-attainment. A sequence of strongly convex regularized problems $(P_k)$ with objectives $f_k(\bx) = x_1 + x_2 + \frac{1}{2k} \|\bx\|^2$ is solved. The plot shows the value of the original objective $x_1 + x_2$ at the solution of each regularized problem. As the regularization parameter $k$ increases (i.e., the regularization term vanishes), the optimal value of the regularized problem converges to the infimum $0$ of the original problem $(P)$, confirming that the value is finite but unattainable.}
    \label{fig:convergence}
\end{figure}

The convergence of the objective value to zero is numerically verified in \Cref{fig:convergence}. We solve a sequence of regularized versions of $(P)$. The solutions of these well-posed problems form a sequence whose objective values approach the infimum, providing clear computational evidence of non-attainment in the original problem.

\subsubsection{Reformulation with Auxiliary Variables}

Sometimes, the problem can be reformulated into an equivalent form that is solvable. SOCPs are often reformulated using auxiliary variables to represent the norms. For our example, the constraint $x_1 \geq \sqrt{1 + x_2^2}$ is equivalent to the existence of $t_1, t_2 \in \R$ such that:
\[
x_1 \geq t_1, \quad t_1 \geq \sqrt{1 + t_2^2}, \quad t_2 = x_2.
\]
This reformulation does not directly help, as the non-polyhedral constraint remains. However, it might open up possibilities for further approximation or decomposition.

\subsubsection{Goal Programming and $\epsilon$-Solutions}

If attainment is not critical and a sufficiently good approximate solution is acceptable, one can simply accept the $\epsilon$-optimal solutions provided by the solver. The solver will terminate when it finds a point $\bx_\epsilon$ such that $f(\bx_\epsilon) \leq \val(P) + \epsilon$. For our example, such points exist for any $\epsilon > 0$. The modeler must then be aware that the solution is not exact but approximate.

\subsection{Conclusion for the Practitioner}

The key takeaway for practitioners is that \emph{robust optimization models with ellipsoidal uncertainty are not guaranteed to be solvable even when they are feasible and bounded}. Care should be taken to check for this possibility, especially in problems where the uncertainty set is large or the constraints are tight. When in doubt, employing a polyhedral approximation is a robust and tractable way to ensure the problem is well-posed and can be solved to true optimality with standard LP solvers.

\section{Conclusion}\label{sec:conclusion}

This paper presented a counterexample in second-order cone programming that challenges the direct extension of a classical linear programming solvability theorem to the conic setting. We constructed a simple SOCP, derived from a robust linear optimization problem with ellipsoidal uncertainty, that is feasible and has a finite optimal value but possesses no optimal solution. This pathology occurs despite the objective function being copositive on the recession cone of the feasible set—a condition that would guarantee solvability for linear programs over polyhedral sets.

Through a detailed analysis, we traced the root cause of this non-attainment to the non-polyhedral nature of the second-order cone constraint, which leads to a feasible set whose image under the linear objective function is not closed. This results in an \emph{attainment gap} where the infimum is not achieved. We further interpreted the example explicitly within the framework of robust optimization, highlighting its practical relevance for modelers who employ ellipsoidal uncertainty sets.

The implications of this work are both theoretical and practical. Theoretically, it underscores the importance of constraint qualifications and the delicate interplay between geometry and optimization in non-polyhedral settings. Practically, it serves as a cautionary note for users of robust optimization, warning them that standard intuitions from linear programming may not hold and that their models might be ill-posed in subtle ways.

To mitigate these issues, we discussed several strategies, including polyhedral approximation and regularization, which can recover solvability and ensure reliable computation. We hope that this analysis fosters a more nuanced understanding of conic programs in robust optimization and provides useful guidance for the formulation and solution of these important models.

Future work could explore the prevalence of this phenomenon in larger, more complex robust optimization problems and investigate specialized algorithms for detecting and handling non-attainment in conic programming.

\section*{Statements and Declarations}
The author received no financial support for the research, authorship, and/or publication of this article. The author declares that there is no conflict of interest regarding the publication of this paper.
\bibliographystyle{plain}

\end{document}